\documentclass[12pt]{article}
\usepackage{mathrsfs}
\usepackage{stmaryrd}
\usepackage{amsfonts,amsmath,amssymb,amscd}
\usepackage{shadow}
\usepackage{graphicx}
\usepackage{color}
\usepackage{longtable}
\usepackage{pstricks,multido}
\usepackage{hyperref}
\usepackage{qtree}
\allowdisplaybreaks

\newtheorem{definition}{Definition }
\newtheorem{mytheorem}{Theorem}

\newtheorem{lemma}[mytheorem]{Lemma}

\newtheorem{remark}{Remark}
\newtheorem{proposition}{Proposition}

\numberwithin{equation}{section}
\numberwithin{figure}{section}

\begin{document}
\begin{center}
{\large{The obstacle problem for nonlinear degenerate
equations with $L^{1}$-data}}
 \end{center}

\begin{center}
{\small Jun Zheng$^1$, Binhua Feng$^2$, Zhihua Zhang$^3$
}

\vskip 3mm
$^{1}$Department of Basic Course,
 Southwest Jiaotong University\\
  Emeishan 614202, China\\
  $^{2}$College of Mathematics and Statistics, Northwest Normal University\\
   Lanzhou 730070, China\\
   $^{3}$School of Mathematical Sciences\\ University of Electronic Science and Technology of China\\
    Chengdu 611731, China
\\[3mm]

\vskip 4mm

$^1$zhengjun2014@aliyun.com, $^2$binhuaf@163.com,
$^3$zhihuamath@aliyun.com
\end{center}

\begin{abstract}
The aim of this paper is to study the obstacle problem with an elliptic operator having degenerate coercivity. We prove the existence of
an entropy solution to the obstacle problem under the assumption of $L^{1}-$summability on the data. Meanwhile, we prove that every entropy solution belongs to some Sobolev space $W^{1,q}(\Omega)$.
\end{abstract}



\section{Introduction}
\subsection{Some remarks and comments}
We begin with some remarks about the Dirichlet problem
\begin{align}
\left\{
\begin{array}{l l}
   -\text {div}\ \frac{a(x,\nabla u)}{(1+|u|)^{\theta(p-1)}}+b|u|^{r-2}u=f, &  \mbox{in }\ \Omega,\\
  u=0, &  \mbox{on }  \ \partial\Omega,
  \end{array} \right.
\label{1.1}
\end{align}
where $\Omega\subset \mathbb{R}^{N}(N\geq 2)$ is a bounded domain, $\theta>0,p>1,1\leq r<p$, $b$ is a constant, $f$ is a measurable function, $a:\Omega \times \mathbb{R}^{N} \rightarrow
\mathbb{R}^{N}$ is a Carath\'{e}odory function, satisfying \begin{align*}
&a(x,\xi)\cdot \xi \geq \alpha |\xi|^{p},\\
&|a(x,\xi)|\leq \beta|\xi|^{p-1},\\
&(a(x, \xi)-a(x, \eta))(\xi-\eta)>0,
\end{align*}
for almost every $x$ in $\Omega$ and for every $\xi,\eta$ in $\mathbb{R}^{N}(\xi\neq\eta)$, where $\alpha$ and $\beta$ are positive constants.

 Problems of $p-$Laplacian type with lower order terms, i.e., $\theta=0$, have been well studied in both the existence and regularity aspects with $f$ having different summability, see \cite{B2006,BG1993,C1995} for instance. But due to the lack of coercivity of the differential operator with $\theta>0$ (it may not tend to infinity on the same space as $u$ becomes large), the classical methods used in order to prove the existence of a solution for elliptic equations (see \cite{LL1965}) cannot be applied even if the data $f$ is regular. However, in \cite{BDO1998}(see also \cite{GP2001,GP2003,P1998}), a whole range of existence results
was proved for the problem
\begin{align}
\left\{
\begin{array}{l l}
   -\text {div}\ \frac{A(x)\nabla u}{(1+|u|)^{\theta}}=f, &  \mbox{in }\ \Omega,\\
  u=0, &  \mbox{on }  \ \partial\Omega,
  \end{array} \right.
\label{1.1'}
\end{align}
under the assumptions
that $A$ is a uniformly elliptic bounded matrix and $\theta \in (0,1]$, yielding solutions in some Sobolev space $W^{1,q}_{0}(\Omega)(q\leq  2)$, if $f$ is regular enough.

For $\theta>0$ and different summability of the datum $f$, Alvino, Boccardo, Ferone, Orsina, Trombetti, et al., have done a lot of work on the existence and regularity of a solution to the problems like \eqref{1.1'}(see \cite{ABFOT2003,AFT1998,BDO1998,GP2001,GP2003,P1998} and references therein). The lack of coercivity implies that the standard Leray-Lions surjectivity theorem can not be used even in the case $f\in W^{-1,p'}(\Omega)$. However, by ¡°cutting¡± the nonlinearity and using the technique of approximation, a pseudomonotone and coercive differential operator on $W^{1,p}_{0}(\Omega)$ can be applied to establish a priori estimates on approximating solutions. As a result, existence and regularity of a solution to the problems like \eqref{1.1'} can be obtained by taking limitation, due to the almost everywhere convergence for the gradients of the approximating solutions.

 This frame work has been extended to the problems with a lower order term, considering
 \begin{align}
\left\{
\begin{array}{l l}
   -\text {div}\ \frac{A(x)\nabla u}{(1+|u|)^{\theta}}+u=f, &  \mbox{in }\ \Omega,\\
  u=0, &  \mbox{on }  \ \partial\Omega,
  \end{array} \right.
\label{1.1''}
\end{align}
with $f\in L^{m}(\Omega),m\geq 1$. Existence of solutions (including distributional solutions ($W^{1,1}$-regular) and entropy solutions) was established in \cite{BB2003,BCO2012,DD2010}. Recently, problems like \eqref{1.1''} were extended to variable Sobolev space in \cite{LG2015}, obtaining renormalized and entropy solutions with $f\in L^{1}(\Omega)$.

However, there are little literatures that consider much higher regularity for entropy solutions of the problems like \eqref{1.1} with noncoercivity, lower terms and $L^{1}-$data. It is the purpose in this paper to establish a certain regularity for these entropy solutions. The interesting cases are those of $1<p<N$, since for $p>N$ the variational
methods of Leray-Lions (see, for instance, \cite{L1969}) can be easily applied to get a solution in some Sobolev space $W^{1,q}(\Omega)$.

  We should mention that problems like \eqref{1.1} with $\theta=0$ and $b=0$, i.e., the case of $p$-Laplacian type, have been well studied in the past and have been extended to unilateral problems in \cite{BC1999,BG1990,P2004}
etc. in Sobolev space, and in \cite{CLR2012,RSU2008} in variable Sobolev space and Orlicz Sobolev space respectively. We should also note that some regularizing effect on the solutions by lower order terms have been studied in \cite{BCO2012,C2010,DD2010} etc..

\subsection{Our problem and main result}
 In this paper, we study the obstacle
problem for nonlinear noncoercive
elliptic equations with lower order term and data $f\in L^{1}(\Omega)$, considering the operator
\begin{eqnarray}
Au=-\text {div}\ \frac{a(x, \nabla u)}{(1+|u|)^{\theta(p-1)}}+b|u|^{r-2}u,\label{Operator A}
\end{eqnarray}
where $b$ is a nonnegative constant and $a:\Omega \times \mathbb{R}^{N} \rightarrow
\mathbb{R}^{N}$ is a Carath\'{e}odory function, satisfying the following conditions:
\begin{align}
&a(x,\xi)\cdot \xi \geq \alpha |\xi|^p,\ \label{2}\\
&|a(x,\xi)|\leq \beta(j(x)+|\xi|^{p-1}),\label{3}\\
&(a(x, \xi)-a(x, \eta))(\xi-\eta)>0,
  \label{5}
\end{align}
for almost every $x$ in $\Omega$, and for every $\xi,\eta \ $ in $\mathbb{R}^{N}$ with $\xi\neq \eta$, where $\alpha,\beta >0$ are constants, and $j$ is a nonnegative function in $ L^{p'}(\Omega)$.

 Given functions
$g,\ \psi \in W^{1,p}(\Omega)$, the obstacle problem associated with \eqref{Operator A} can be formulated in terms of the inequality
 \begin{align}
\int_{\Omega}\frac{a(x, \nabla u)}{(1+|u|)^{\theta(p-1)}}\cdot\nabla (u-v)\text{d}x&+\int_{\Omega}b|u|^{r-2}u(u-v)\text{d}x\notag\\&\leq \int_{\Omega}f(u-v)\text{d}x,\ \ \forall v\in K_{g,\psi},\label{classical obstacle}
\end{align}
whenever $1\leq r< p,f\in W^{-1,p'}(\Omega)$ and the convex subset
\begin{align*}
K_{g,\psi}={ \{v\in W^{1,p}(\Omega);\ v-g\ \in W^{1,p}_{0}(\Omega),\
v\geq \psi,\ \text{a.e. in}\ \Omega\} }
\end{align*}
is nonempty. However, for a general $f\in L^{1}(\Omega)$, the third integration in \eqref{classical obstacle} is not
defined, and following \cite{ABFOT2003,BBGGPV1995,BC1999} etc., we are led to the more general definition of a solution to
the obstacle problem, using the truncation function
\begin{align*}
T_{s}(t)=\max\{-s,\min\{s,t\}\},\ \ s,t\in \mathbb{R}.
\end{align*}

\begin{definition}\label{Definition 1}
An entropy solution of the obstacle problem associated with $(f,\psi,g)$ and $f\in L^{1}(\Omega)$ is a measurable function $u$ such that $u\geq \psi$ a.e. in $\Omega$, $|u|^{r-1}\in L^{1}(\Omega)$ and for every $s>0$, $T_{s}(u)-T_{s}(g)\in W_{0}^{1,p}(\Omega)$ and
\begin{align}\label{1.9}
\int_{\Omega}\!\!\frac{a(x, \nabla u)}{(1+|u|)^{\theta(p-1)}}\cdot\nabla (T_{s}(u-v))\text{d}x+\int_{\Omega}b|u|^{r-2}uT_s(u-v)\text{d}x \notag\\
\leq \int_{\Omega}fT_{s}(u-v)\text{d}x,\ \ \forall v\in K_{g,\psi}.
\end{align}
\end{definition}

  Throughout this paper, without special statements, we always assume
 \begin{align*}
 2-\frac{1}{N}<p<N,\ 1\leq r<p,\ 0\leq \theta< \min\bigg\{\frac{N}{N-1}-\frac{1}{p-1},\frac{p-r}{p-1}\bigg\},
 \end{align*}
and
 $\psi,g\in W^{1,p}(\Omega)\cap L^{\infty}(\Omega)$, and $ (\psi-g)^{+}\in W^{1,p}_0(\Omega)$ such that $K_{g,\psi}\neq \emptyset$. We make further assumption on $a$:
\begin{align}
|a(x, \xi)-a(x, \eta)|\leq \gamma
\left\{
\begin{array}{l l}
  |\xi-\eta|^{p-1}, &  \mbox{if $1<p< 2$,}\\
  (1+|\xi|+|\eta|)^{p-2}|\xi-\eta|, &  \mbox{if $p\geq 2$,}
  \end{array} \right.
\label{6}
\end{align}
for almost every $x$ in $\Omega$, and for every $\xi,\eta$ in $\mathbb{R}^{N}$, where $\gamma $ is a positive constant.

 The main result in this paper is
\begin{mytheorem}\label{Theorem 1.1}
(i)\ Let $f\in L^{1}(\Omega)$. Then there exists at least one entropy solution $u$ of the obstacle problem associated with $(f,\psi,g) $. In addition, $u$ depends continuously on $f$, i.e., if $f_{n}\rightarrow f$ in $L^{1}(\Omega)$ and $u_n$ is a solution to the obstacle problem associated with $(f_n,\psi,g)$, then
\begin{align}\label{Q}
u_n\rightarrow u\ \text{in}\ W^{1,q}(\Omega),\forall q
\in\! \left\{
\begin{array}{l l}
  \!\!\!\big( \frac{N(r-1)}{N+r-1}, \frac{N(p-1)(1-\theta)}{N-1-\theta(p-1)}\big), &  \!\mbox{\!if}\ \  \frac{2N-1}{N-1}\leq r<p,\\
  \!\!\!(1,\frac{N(p-1)(1-\theta)}{N-1-\theta(p-1)} ), &\!\!\!\!\!\!\!\!\!\! \text{\!if}\  \ 1\leq \!r\!< \min\{\frac{2N-1}{N-1},p\}.
  \end{array} \right.
\end{align}
(ii) Let $f\in L^{1}(\Omega)$. Then every entropy solution $u$ of the obstacle problem associated with $(f,\psi,g) $ belongs to Sobolev space $W^{1,q}(\Omega)$ for all $ q$ given by \eqref{Q}.
  \end{mytheorem}

\begin{remark}\label{Remark 1}
 (i) \ \ $\bigg( \frac{N(r-1)}{N+r-1}, \frac{N(p-1)(1-\theta)}{N-1-\theta(p-1)}\bigg)\subset \bigg(1,\frac{N(p-1)(1-\theta)}{N-1-\theta(p-1)} \bigg)$ if $ \frac{2N-1}{N-1}\leq r<p$.
Indeed, $ \theta<\frac{p-r}{p-1}+\frac{p(r-1)}{N(p-1)}\Leftrightarrow \frac{N(p-1)(1-\theta)}{N-1-\theta(p-1)}>\frac{N(r-1)}{N+r-1}$, while $ \frac{2N-1}{N-1}\leq r $ gives $\frac{N(r-1)}{N+r-1}\geq 1 $. Thus $u_n\rightarrow u\ \text{in}\ W^{1,q}(\Omega)\ \ for\ all \ q\in \big(1,\frac{N(p-1)(1-\theta)}{N-1-\theta(p-1)}\big)$.

 (ii)\ \ $r-1<\frac{Nq}{N-q}$. Indeed, by $ 1\leq r<\frac{2N-1}{N-1}$, there holds $r-1<\frac{N}{N-1}< \frac{Nq}{N-q}$ for any $q>1$, particularly, for $ q\in (1,\frac{N(p-1)(1-\theta)}{N-1-\theta(p-1)} )$. For $r\geq \frac{2N-1}{N-1}$, it suffices to note that $ q> \frac{N(r-1)}{N+r-1}\Leftrightarrow r-1< \frac{Nq}{N-q}$.

 (iii)\ \  $q<p$. Indeed, $0\leq \theta<\frac{N}{N-1}-\frac{1}{p-1}<\frac{N-1}{p-1} \Rightarrow\frac{N(p-1)(1-\theta)}{N-1-\theta(p-1)}<p $.
\end{remark}

\subsection{Notations}
  ~~~~$\|u\|_{p}=\|u\|_{L^{p}(\Omega)}$ is the norm of $L^{p}(\Omega)$, where $1\leq p\leq \infty$.

  $\|u\|_{1,p}=\|u\|_{W^{1,p}(\Omega)}$ is the norm of $W^{1,p}(\Omega)$, where $1<p<\infty$.

  $p'=\frac{p}{p-1}$ with $1<p<\infty$.

            $u^+=\max\{u,0\}$, $u^-=(-u)^+$, if $u$ is a real-valued function.

            $C$ is a constant, which may be different from each other.

            $\{u>s\}=\{x\in \Omega;u(x)>s\}$.

            $\{u\leq s\}=\Omega\setminus\{u>s\}$.

   $\{u<s\}=\{x\in \Omega;u(x)<s\}$.

    $\{u\geq s\}=\Omega\setminus\{u<s\}$.

      $\{u=s\}=\{x\in \Omega;u(x)=s\}$.

      $\{t\leq u < s\}=\{u\geq t\}\cap \{u<s\} $.

       $\mathcal {L}^N$ is the Lebesgue measure of $\mathbb{R}^N$.

       $|E|=\mathcal {L}^N(E)$ for a measurable set $E$.

\section{Lemmas on entropy solutions}
It is worthy to note that, for any smooth function $f_n$, there exists at least one solution to the obstacle problem \eqref{classical obstacle}. Indeed, one can proceed exactly as Theorem 7.1 of \cite{DD2010} to obtain $W^{1,p}-$solutions due to the assumptions \eqref{2}-\eqref{5} on $a$ and $r-1<p$. These solutions, in particular, are also entropy solutions. In this section we establish several auxiliary results on convergence of sequences
of entropy solutions when $f_{n}\rightarrow f\ \text{in}\ L^1(\Omega)$. We start with an a priori estimate.

\begin{lemma}\label{Lemma 3.1}
Let $v_{0}\in K_{g,\psi}\cap L^{\infty}(\Omega)$, and let $u$ be an entropy solution of the obstacle problem associated with $(f,\psi,g)$. Then, we have
\begin{align*}
\int_{\{|u|<t\}}\!\!\frac{|\nabla u|^{p}}{(1+|u|)^{\theta(p-1)}}\text{d}x\leq C(1+t^r),\ \forall t\!>\!0,
\end{align*}
where $C$ is a positive constant depending only on $\alpha,\beta,p,r,b,\|j\|_{p'},\|\nabla v_{0}\|_{p},\|v_{0}\|_{\infty},$ and $\|f\|_{1}$.
\end{lemma}
\textbf{Proof.}\ \ Take $v_{0}$ as a test function in \eqref{1.9}. For $t$ large enough such that $t-\|v_0\|_{\infty}>0$, we get
\begin{align}
\int_{\{|v_{0}-u|<t\}}\frac{a(x, \nabla u)\cdot\nabla u}{(1+|u|)^{\theta(p-1)}}\text{d}x\leq \int_{\{|v_{0}-u|<t\}}\frac{a(x, \nabla u)\cdot\nabla v_{0}}{(1+|u|)^{\theta(p-1)}}\text{d}x
\notag\\
   +\int_{\Omega}f(T_{t}(u-v_0))\text{d}x-\int_{\Omega}b|u|^{r-2}uT_{t}(u-v_0))\text{d}x.\label{1.90}
\end{align}
We estimate each integration in the right-hand side of \eqref{1.90}. It follows from \eqref{3} and Young's inequality with $\varepsilon>0$ that
\begin{align}
&\int_{\{|v_{0}-u|<t\}}\frac{a(x, \nabla u)\cdot\nabla v_{0}}{(1+|u|)^{\theta(p-1)}}\text{d}x\leq \int_{\{|v_{0}-u|<t\}}\!\!\!\!\!\!\frac{\beta(|j|+|\nabla u|^{p-1})\cdot|\nabla v_{0}|}{(1+|u|)^{\theta(p-1)}}\text{d}x\notag\\
&\leq \int_{\{|v_{0}-u|<t\}}\frac{\beta \varepsilon(|j|^{p'}+|\nabla u|^{p})}{(1+|u|)^{\theta(p-1)}}\text{d}x+\int_{\{|v_{0}-u|<t\}}\frac{\beta C(\varepsilon)|\nabla v_{0}|^{p}}{(1+|u|)^{\theta(p-1)}}\text{d}x\notag\\
&\leq  \varepsilon \int_{\{|v_{0}-u|<t\}}\frac{ |\nabla u|^p}{(1+|u|)^{\theta(p-1)}}\text{d}x+C(\|j\|_{p'}^{p'}+\|\nabla v_{0}\|_{p}^{p}).\label{1.91}
\end{align}
\begin{align}
-\int_{\Omega}b|u|^{r-2}uT_{t}(u-v_0)\text{d}x&=-\int_{\{|u-v_0|\leq t\}} b|u|^{r-2}uT_{t}(u-v_0)\text{d}x\notag\\&
-\int_{\{|u-v_0|> t\}}b |u|^{r-2}uT_{t}(u-v_0)\text{d}x. \label{1.92}
\end{align}
Note that on the set $ \{|u-v_0|\leq t\}$,
\begin{align}
 ||u|^{r-2}uT_{t}(u-v_0)|\leq t|t+\|v_{0}\|_{\infty}|^{r-1}\leq C(1+t^{r}),\label{1.93}
\end{align}
where $C$ is a constant depending only on $r,\|v_{0}\|_{\infty}$.\\
On the set $ \{|u-v_0|> t\}$, we have $|u|\geq t-\|v_{0}\|_{\infty} >0$, thus $u$ and $T_{t}(u-v_0)$ have the same sign. It fllows
\begin{align}
-\int_{\{|u-v_0|> t\}} b|u|^{r-2}uT_{t}(u-v_0)\text{d}x\leq 0.\label{1.94}
\end{align}
Combining \eqref{1.92}-\eqref{1.94}, we get
\begin{align}
-\int_{\Omega}b|u|^{r-2}uT_{t}(u-v_0)\text{d}x\leq C(1+t^{r}). \label{1.95}
\end{align}
\begin{align}
\int_{\{|v_{0}-u|<t\}}\frac{|\nabla u|^p}{(1+|u|)^{\theta(p-1)}}\text{d}x&\leq C(\|j\|_{p'}^{p'}+\|\nabla v_{0}\|_{p}^{p}+t\|f\|_{1}+1+t^{r})\notag\\
&\leq C(1+t^{r}).\label{3.3}
\end{align}
Replacing $t$ with $t+\|v_0\|_{\infty}$ in \eqref{3.3} and noting that $\{|u|<t\}\subset \{|v_0-u|<t+\|v_0\|_{\infty}\}$, one may obtain the desired result.\ $\square$\\

 In the rest of this section, let $\{u_{n}\}$ be a sequence of entropy solutions of the obstacle problem associated with $(f_{n},\psi,g)$ and assume that
\begin{align*}
f_{n}\rightarrow f\ \ \text{in}\ L^1(\Omega) \ \ \text{and}\ \ \|f_{n}\|_1\leq \|f\|_{1}+1.
\end{align*}

\begin{lemma}\label{Lemma 3.2}
There exists a measurable function $u$ such that
$
u_{n}\rightarrow u
$  in  measure, and $
T_{k}(u_n)\rightharpoonup T_{k}(u)$  weakly in $W^{1,p}(\Omega)$ for any $ k>0$. Thus
$ T_{k}(u_n)\rightarrow T_{k}(u) $ strongly in $L^{p}(\Omega)$ and\ a.e. in $\Omega$.
\end{lemma}

\textbf{Proof.}\ \ Let $s,t$ and $\varepsilon$ be positive numbers. One may verify that for every $m,n\geq 1$,
\begin{align}
\mathcal {L}^{N}(\{|u_{n}-u_{m}|>s\})&\leq \mathcal {L}^{N}(\{|u_{n}|>t\})+
\mathcal {L}^{N}(\{|u_{m}|>t\})\notag\\
&+\mathcal {L}^{N}(\{|T_{k}(u_n)-T_{k}(u_m)|>s\}),\label{3.5}
\end{align}
and
\begin{align}
\mathcal {L}^{N}(\{|u_{n}|>t\})=\frac{1}{t^p}\int_{\{|u_n|>t\}}t^p\text{d}x\leq
\frac{1}{t^p}\int_{\Omega}|T_t(u_n)|^p\text{d}x.
\label{3.6}
\end{align}
Due to $v_{0}=g+(\psi-g)^{+}\in K_{g,\psi}\cap L^{\infty}(\Omega)$, by Lemma~\ref{Lemma 3.1}, one has
\begin{align}
\int_{\Omega}\!|\nabla T_t(u_n)|^p\text{d}x&\!=\!\int_{\{|u_n|<t\}}|\nabla u_n|^p\text{d}x\notag\\
&\leq \!C(1\!+\!t)^{\theta(p-1)}(1+t^{r}).\label{3.7}
\end{align}
Note that $ T_{t}(u_{n})-T_{t}(g)\in W^{1,p}_{0}(\Omega)$. by \eqref{3.6}, \eqref{3.7} and Poincar\'{e} inequality, for every $t>\|g\|_{\infty}$ and for some positive constant $C$ independent of $n$ and $t$, there holds
\begin{align*}
\mathcal {L}^{N}(\{|u_{n}|>t\})&\leq \frac{1}{t^p}\int_{\Omega}|T_t(u_n)|^p\text{d}x\notag\\
&\leq \frac{2^{p-1}}{t^p}\int_{\Omega}|T_t(u_n)-T_{t}(g)|^p\text{d}x+\frac{2^{p-1}}{t^p}\|g\|_{p}^{p}\notag\\
&\leq \frac{C}{t^p}\int_{\Omega}|\nabla T_t(u_n)-\nabla T_{t}(g)|^p\text{d}x+\frac{2^{p-1}}{t^p}\|g\|_{p}^{p}\notag\\
&\leq\frac{C}{t^p}\int_{\Omega}|\nabla T_t(u_n)|^p\text{d}x+\frac{C}{t^p}\|g\|_{1,p}^{p}\notag\\
&\leq \frac{C(1+t^{r+\theta(p-1)})}{t^p}.
\end{align*}
Since $0\leq \theta<\frac{p-r}{p-1}$, there exists $t_{\varepsilon}>0$ such that
\begin{align}
\mathcal {L}^{N}(\{|u_{n}|>t\})< \frac{\varepsilon}{3},\ \ \forall\ t\geq t_{\varepsilon},\ \forall\ n\geq 1.\label{3.8}
\end{align}
Now we have as in \eqref{3.6}
\begin{align}
\mathcal {L}^{N}(\{|T_{t_{\varepsilon}}(u_{n})-T_{t_{\varepsilon}}(u_{m})|>s\})
&=\frac{1}{s^p}\int_{\{|T_{t_{\varepsilon}}
(u_{n})-T_{t_{\varepsilon}}(u_{m})|>s\}}
s^p\text{d}x\notag\\
&\leq\frac{1}{s^p}\int_{\Omega}
|T_{t_{\varepsilon}}(u_{n})-T_{t_{\varepsilon}}(u_{m})|^p\text{d}x.\label{3.9}
\end{align}
Using \eqref{3.7} and the fact that $T_{t}(u_{n})-T_{t}(g)\in W^{1,p}_{0}(\Omega)$ again, we see that $\{T_{t_{\varepsilon}}(u_{n})\}$ is a bounded sequence in $W^{1,p}(\Omega)$. Thus, up to a subsequence, $ \{T_{t_{\varepsilon}}(u_{n})\}$ converges strongly in $L^{p}(\Omega)$. Taking into account \eqref{3.9}, there exists
$n_{0}=n_{0}(t_{\varepsilon},s)\geq 1$ such that
\begin{align}
\mathcal {L}^{N}(\{|T_{t_{\varepsilon}}(u_{n})-T_{t_{\varepsilon}}(u_{m})|>s\})< \frac{\varepsilon}{3},\ \ \forall\ n,m\geq n_{0}.\label{3.10}
\end{align}
Combining \eqref{3.5}, \eqref{3.8} and \eqref{3.10}, we obtain
\begin{align*}
\mathcal {L}^{N}(\{|u_{n}-u_{m}|>s\})< \varepsilon,\ \ \forall\ n,m\geq n_{0}.
\end{align*}
Hence $\{u_{n}\}$ is a Cauchy sequence in measure, and therefore there exists a measurable function $u$ such that $u_n\rightarrow u$ in measure. The remainder of the lemma is a consequence of the fact
that  $\{T_{k}(u_{n})\}$ is a bounded sequence in $W^{1,p}(\Omega)$. \ $\square$

\begin{proposition}\label{Proposition 3.1}
There exists a subsequence of $\{u_n\}$ and a measurable function $u$ such that for each $q$ given in \eqref{Q}, we have
 \begin{align*}
u_n\rightarrow u \ \ \text{strongly\ in} \ W^{1,q}(\Omega).
\end{align*}
If moreover $0\leq \theta< \min\{\frac{1}{N-p+1},\frac{N}{N-1}-\frac{1}{p-1},\frac{p-r}{p-1}\}$, then
 \begin{align*}
\frac{a(x, \nabla u_n)}{(1+|u_n|)^{\theta(p-1)}} \rightarrow \frac{a(x, \nabla u)}{(1+|u|)^{\theta(p-1)}} \ \ \text{strongly\ in} \ (L^1(\Omega))^{N}.
\end{align*}
\end{proposition}

 To prove Proposition \ref{Proposition 3.1}, we need two preliminary lemmas.
\begin{lemma}\label{Lemma 3.4}
There exists a subsequence of $\{u_n\}$ and a measurable function $u$ such that for each $q$ given in \eqref{Q}, we have
 $
u_n\rightharpoonup u \ \text{weakly\ in} \ W^{1,q}(\Omega),$ and $
u_n\rightarrow u \ \text{strongly\ in} \ L^q(\Omega).
$
\end{lemma}
\textbf{Proof.}\ \ Let $k>0$ and $n\geq 1$. Define $D_{k}=\{|u_{n}|\leq k\}$ and $B_{k}=\{k\leq|u_{n}|< k+1\}$. Using Lemma \ref{Lemma 3.1} with $v_{0}=g+(\psi-g)^{+}$, we get
\begin{align}
\int_{D_{k}}\frac{|\nabla u_{n}|^p}{(1+|u_n|)^{\theta(p-1)}}\text{d}x\leq C(1+k^r),\label{3.11}
\end{align}
where $C$ is a positive constant depending only on $\alpha,\ \beta,\ p ,\ r, \|j\|_{p'},\ \|f\|_{1},\ \|\nabla v_{0}\|_{p},$ and $ \|v_{0}\|_{\infty}$.

 Using the function $T_{k}(u_{n})$ for $k> \|g\|_{\infty},\|\psi\|_{\infty}$, as a test function for the problem associated with $(f_{n},\psi,g)$, we obtain
\begin{align*}
\int_{\Omega}\frac{a(x,\nabla u_n)\cdot \nabla (T_{1}(u_n-T_{k}(u_n)))}{(1+|u_n|)^{\theta(p-1)}}\text{d}x+
\int_{\Omega}b|u_n|^{r-2}u_nT_{1}(u_n-T_{k}(u_n))\text{d}x\notag\\\leq \int_{\Omega}f_nT_{1}(u_n-T_{k}(u_n))\text{d}x,
\end{align*}
which and \eqref{2} give
\begin{align*}
\int_{B_{k}}\frac{\alpha|\nabla u_{n}|^p}{(1+|u_n|)^{\theta(p-1)}}\text{d}x+\int_{\Omega}b|u_n|^{r-2}
u_nT_{1}(u_n-T_{k}(u_n))\text{d}x&\leq \|f_n\|_{1}\notag\\&\leq \|f\|_{1}+1.
\end{align*}
Note that on the set $\{|u_n|\geq k+1\} $, $u_n$ and $T_{1}(u_n-T_{k}(u_n))$ have the same sign. Then
\begin{align*}
\int_{\Omega}|u_n|^{r-2}u_nT_{1}(u_n-T_{k}(u_n))\text{d}x
&=\int_{D_{k}}|u_n|^{r-2}u_nT_{1}(u_n-T_{k}(u_n))\text{d}x\notag\\
&+\int_{B_{k}}|u_n|^{r-2}u_nT_{1}(u_n-T_{k}(u_n))\text{d}x
\notag\\
&+\int_{\{|u_n|\geq k+1\}}|u_n|^{r-2}u_nT_{1}(u_n-T_{k}(u_n))\text{d}x
\notag\\
&\geq \int_{B_{k}}|u_n|^{r-2}u_nT_{1}(u_n-T_{k}(u_n))\text{d}x.
\end{align*}
Thus we have
\begin{align}
\int_{B_{k}}\frac{\alpha|\nabla u_{n}|^p}{(1+|u_n|)^{\theta(p-1)}}\text{d}x+&\leq \|f\|_{1}+1-\int_{B_{k}}b|u_n|^{r-2}u_nT_{1}(u_n-T_{k}(u_n))\text{d}x\notag\\
&\leq \|f\|_{1}+1+\int_{B_{k}}b|u_n|^{r-1}\text{d}x\notag\\
&\leq C\bigg(1+\bigg(\int_{B_{k}}|u_n|^{q^{*}}\text{d}x
\bigg)^{\frac{r-1}{q^{*}}}|B_{k}|^{1-\frac{r-1}{q^{*}}}\bigg)
,\label{3.12}
\end{align}
where $q$ is given in \eqref{Q} and $q^*=\frac{Nq}{N-q} $.\\ Let $s=\frac{q\theta(p-1)}{p}$. Note that $q<p$ and $\frac{ps}{p-q}<q^*$. For $\forall k>0$, we estimate $\int_{B_{k}}|\nabla u_n|^q\text{d}x$ as follows
\begin{align*}
\int_{B_{k}}\!\!|\nabla u_n|^q\text{d}x&\!\!=\!\!\int_{B_{k}}\frac{|\nabla u_{n}|^q}{(1+|u_n|)^{s}}\cdot (1+|u_n|)^{s}\text{d}x\notag\\
&\!\!\leq \bigg(\int_{B_{k}}\frac{|\nabla u_{n}|^p}{(1+|u_n|)^{\theta(p-1)}}\text{d}x\bigg)^{\frac{q}{p}}
\bigg(\int_{B_{k}}(1+|u_n|)^{\frac{ps}{p-q}}\text{d}x\bigg)^{\frac{p-q}{p}}\notag\\
&\!\!\leq \!C\bigg(\!\!\int_{B_{k}}\!\!\frac{|\nabla u_{n}|^p}{(1+|u_n|)^{\theta(p-1)}}\text{d}x\bigg)^{\!\!\frac{q}{p}}
\bigg(|B_{k}|^{\frac{p-q}{p}}+
\bigg(\!\!\int_{B_{k}}\!\!|u_n|^{\frac{ps}{p-q}}\text{d}x\bigg)^{\!\!\frac{p-q}{p}}\bigg)\notag\\
&\!\!\leq \!C\bigg(\!\!\int_{B_{k}}\!\!\frac{|\nabla u_{n}|^p}{(1+|u_n|)^{\theta(p-1)}}\text{d}x\!\bigg)^{\!\frac{q}{p}}\!\!
\bigg(\!|B_{k}|^{\frac{p-q}{p}}\! +\!\!\bigg(\!\!\int_{B_{k}}\!\!\!|u_n|^{q^*}\!\text{d}x\bigg)^{\!\frac{s}{q^*}}
\!\!|B_{k}|^{\frac{p-q}{p}-\frac{s}{q^*}}\!\bigg)
\notag\\
&\!\!\leq \!C\bigg(\!\!|B_{k}|^{\frac{p-q}{p}}\!\!+|B_{k}|^{\frac{p-q}{p}-\frac{s}{q^*}}
\bigg(\!\!\int_{B_{k}}\!\!|u_n|^{q^*}\text{d}x\!\bigg)^{\frac{s}{q^*}}\!\!\!\!\!\!+|B_{k}|^{1-p_1}\bigg(\!\!\int_{B_{k}}\!\!|u_n|^{q^*}\text{d}x\!\bigg)^{p_1}\notag\\
&+|B_{k}|^{1-p_2}\bigg(\int_{B_{k}}|u_n|^{q^*}\text{d}x\bigg)^{p_2}\bigg)\ \ \ \ \text{by}\  \eqref{3.12}\notag\\
&\!\!= C\bigg(|B_{k}|^{\frac{p-q}{p}}+
|B_{k}|^{\frac{p-q}{p}-\frac{s}{q^*}}\bigg(\int_{B_{k}}|u_n|^{q^*}\text{d}x\bigg)^{\frac{s}{q^*}}\notag\\
&+|B_{k}|^{1-p_1-C_1}|B_{k}|^{C_1}\bigg(\int_{B_{k}}|u_n|^{q^*}\text{d}x\bigg)^{p_1}\notag\\
&+|B_{k}|^{1-p_2-C_2}|B_{k}|^{C_2}\bigg(\int_{B_{k}}|u_n|^{q^*}\text{d}x\bigg)^{p_2}\bigg),
\end{align*}
where $p_1=\frac{q}{p}\frac{r-1}{q^*}, p_2=\frac{s}{q^*}+\frac{q}{p}\frac{r-1}{q^*}$, $C_1$ and $C_2$ are positive constants to be chosen later.\\
Note that $\theta<\frac{p-r}{p-1}$, it follows
\begin{align*}
\frac{\theta(p-1)}{p}+\frac{r-1}{p}<\frac{p-1}{p}<1-\frac{1}{N}=1-\frac{1}{q}+\frac{1}{q^*}.
\end{align*}
Thus
\begin{align*}
&\frac{\theta q(p-1)}{p}+\frac{q(r-1)}{p}+1<q+\frac{q}{q^*}\notag\\
&\Leftrightarrow s+\frac{q(r-1)}{p}+1<q+\frac{q}{q^*}
\notag\\
&\Leftrightarrow p_2+\frac{1-p_2}{q^*+1}<\frac{q}{q^*}.
\end{align*}
Note that $p_1<p_2<1$. Then for $ i=1,2$, we always have
\begin{align*}
p_i+\frac{1-p_i}{q^*+1}<\frac{q}{q^*}<1.
\end{align*}
From this, we may find positive $C_{i}(i=1,2)$ such that
\begin{align}
p_i+\frac{1-p_i}{q^*+1}<p_i+C_{i}<\frac{q}{q^*}<1,\ \ i=1,2. \label{pici}
\end{align}
It follows
\begin{align*}
\frac{1-p_i}{q^*+1}<C_{i}\Leftrightarrow 1-p_i-C_i<C_iq^*,\ \ i=1,2,
\end{align*}
which implies
\begin{align}
C_{i}\alpha_iq^*=\frac{C_iq^*}{1-p_i-C_i}>1,\ \ i=1,2,\label{alpha}
\end{align}
with $\alpha_i=\frac{1}{1-p_i-C_i}>1,\ i=1,2$. Let $\beta_i=\frac{1}{p_i+C_i}>1, i=1,2$. Then we have $\frac{1}{\alpha_i}+\frac{1}{\beta_i}=1(i=1,2)$.\\
Since $|B_{k}|\leq \frac{1}{k^{q^*}}\int_{B_{k}}|u_n|^{q^*}\text{d}x $,  $|B_{k}|^{1-p_1-C_1} \leq |\Omega|^{1-p_1-C_1}$ and $|B_{k}|^{1-p_2-C_2} \leq |\Omega|^{1-p_2-C_2}$, we have for $k\geq k_0\geq 1$
 \begin{align*}
\int_{B_{k}}|\nabla u_n|^q\text{d}x&\leq \frac{C}{k^{q^*\big(\frac{p-q}{p}-\frac{s}{q^*}\big)}}
\bigg(\int_{B_{k}}|u_n|^{q^*}\text{d}x\bigg)^{\frac{p-q}{p}}
+\frac{C}{k^{q^*C_1}}\bigg(\int_{B_{k}}|u_n|^{q^*}\text{d}x\bigg)^{p_1+C_1}\notag\\
&+\frac{C}{k^{q^*C_2}}\bigg(\int_{B_{k}}|u_n|^{q^*}\text{d}x\bigg)^{p_2+C_2}.
\end{align*}
Summing up from $k=k_0$ to $k=K$ and using H\"{o}lder's inequality, one has
 \begin{align}
\sum_{k=k_{0}}^{K}\int_{B_{k}}|\nabla u_n|^q\text{d}x&\leq C
\bigg(\sum_{k=k_{0}}^{K}
\frac{1}{k^{q^*(\frac{p-q}{p}-\frac{s}{q^*})\frac{p}{q}}}
\bigg)^{\frac{q}{p}}\cdot
\bigg(\sum_{k=k_{0}}^{K}\int_{B_{k}}|u_n|^{q^*}
\text{d}x\bigg)^{\frac{p-q}{p}}\notag\\
&+C\bigg(\sum_{k=k_{0}}^{K}
\frac{1}{k^{q^*C_1\alpha_1}}\bigg)^{\frac{1}{\alpha_1}}
\cdot\bigg(\sum_{k=k_{0}}^{K}\bigg(\int_{B_{k}}|u_n|^{q^*}
\text{d}x\bigg)^{\beta_1(p_1+C_1)}\bigg)^{\frac{1}{\beta_1}}\notag\\
&+C\bigg(\sum_{k=k_{0}}^{K}
\frac{1}{k^{q^*C_2\alpha_2}}\bigg)^{\frac{1}{\alpha_2}}
\cdot\bigg(\sum_{k=k_{0}}^{K}\bigg(\int_{B_{k}}|u_n|^{q^*}
\text{d}x\bigg)^{\beta_2(p_2+C_2)}\bigg)^{\frac{1}{\beta_2}}\notag\\
&= C
\bigg(\sum_{k=k_{0}}^{K}
\frac{1}{k^{q^*(\frac{p-q}{p}-\frac{s}{q^*})\frac{p}{q}}}
\bigg)^{\frac{q}{p}}\cdot
\bigg(\sum_{k=k_{0}}^{K}\int_{B_{k}}|u_n|^{q^*}
\text{d}x\bigg)^{\frac{p-q}{p}}\notag\\
&+C\bigg(\sum_{k=k_{0}}^{K}
\frac{1}{k^{q^*C_1\alpha_1}}\bigg)^{\frac{1}{\alpha_1}}
\cdot\bigg(\sum_{k=k_{0}}^{K}\int_{B_{k}}|u_n|^{q^*}
\text{d}x\bigg)^{p_1+C_1}\notag\\
&+C\bigg(\sum_{k=k_{0}}^{K}
\frac{1}{k^{q^*C_2\alpha_2}}\bigg)^{\frac{1}{\alpha_2}}
\cdot\bigg(\sum_{k=k_{0}}^{K}\int_{B_{k}}|u_n|^{q^*}
\text{d}x\bigg)^{p_2+C_2}.
 \label{3.18}
\end{align}
Note that
\begin{align}
\int_{\{|u_n|\leq K\}}|\nabla u_n|^q\text{d}x=\int_{D_{k_0}}|\nabla u_n|^q\text{d}x+\sum_{k=k_{0}}^{K}\int_{B_{k}}|\nabla u_n|^q
\text{d}x.\label{3.19}
\end{align}
To estimate the first integral in the right-hand side of \eqref{3.19}, we compute by using H\"{o}lder's inequality and \eqref{3.11}, obtaining
\begin{align}
\int_{D_{k_0}}\!\!|\nabla u_n|^q\text{d}x
&\leq \bigg(\int_{D_{k_0}}\!\!\frac{|\nabla u_{n}|^p}{(1+|u_n|)^{\theta(p-1)}}\text{d}x\bigg)^{\frac{q}{p}}
\bigg(\int_{D_{k_0}}\!\!(1+|u_n|)^{\frac{ps}{p-q}}
\text{d}x\bigg)^{\frac{p-q}{p}}\notag\\
&\leq C,\label{3.20}
\end{align}
where $C$ depending only on $\alpha,\ \beta,\ p,\ \theta,\  \|j\|_{p'},\ \|f\|_{1},\ \|\nabla v_{0}\|_{p},\ \|v_{0}\|_{\infty}$ and $k_{0}$.

 Note that $ \sum_{k=k_{0}}^{K}
\frac{1}{k^{q^*(\frac{p-q}{p}-\frac{s}{q^*})\frac{p}{q}}}$ and $\sum_{k=k_{0}}^{K}
\frac{1}{k^{q^*C_i\alpha_i}}$ converge due to the fact that $q^*(\frac{p-q}{p}-\frac{s}{q^*})\frac{p}{q}>1$ and $q^*C_i\alpha_i>1$ by \eqref{alpha}, respectively. Combining \eqref{3.18}-\eqref{3.20}, we get for $k_{0}$ large enough
\begin{align}
\int_{\{|u_n|\leq K\}}\!\!\!\!\!\!|\nabla u_n|^q\text{d}x&\leq C+C\bigg(\int_{\{|u_n|\leq K\}}\!\!\!\!|u_n|^{q^*}
\text{d}x\bigg)^{\frac{p-q}{p}}\!\!\!\!\!\!+
C\bigg(\int_{\{|u_n|\leq K\}}\!\!|u_n|^{q^*}
\text{d}x\bigg)^{p_1+C_1}\notag\\
&+C\bigg(\int_{\{|u_n|\leq K\}}|u_n|^{q^*}
\text{d}x\bigg)^{p_2+C_2}.\label{3.21}
\end{align}
Since $p>q$, $T_{K}(u_{n})\in W^{1,q}(\Omega),\ T_{K}(g)=g\in W^{1,q}(\Omega)$ for $K>\|g\|_{\infty}$. Hence $T_{K}(u_{n})-g\in W^{1,q}_0(\Omega)$. Using the Sobolev embedding $W^{1,q}_0(\Omega)\subset L^{q^{*}}(\Omega)$ and Poincar\'{e} inequality, we obtain
\begin{align}
\|T_{K}(u_{n})\|^{q}_{q^*}&\leq 2^{q-1}(\|T_{K}(u_{n})-g\|^{q}_{q^*}+\|g\|^{q}_{q^*})\notag\\
&\leq C(\|\nabla (T_{K}(u_{n})-g)\|^{q}_{q}+\|g\|^{q}_{q^*})\notag\\
&\leq C (\|\nabla T_{K}(u_{n})\|^{q}_{q}+\|\nabla g\|^{q}_{q}+\|g\|^{q}_{q^*})\notag\\
&\leq C\bigg(1+\int_{\{| u_n|\leq K\}}|\nabla u_n|^{q}
\text{d}x\bigg).
\label{3.22}
\end{align}
Using the fact that
\begin{align}
\int_{\{| u_n|\leq K\}}| u_n|^{q^{*}}
\text{d}x\leq \int_{\{| u_n|\leq K\}}|T_{K}( u_n)|^{q^{*}}
\text{d}x\leq \|T_{K}( u_n)\|^{q^*}_{q^*},
\label{3.23}
\end{align}
we obtain from \eqref{3.21}-\eqref{3.22},
\begin{align}
\int_{\{|u_n|\leq K\}}|\nabla u_n|^q\text{d}x&\leq C+C\bigg(1+\int_{\{| u_n|\leq K\}}|\nabla u_n|^{q}
\text{d}x\bigg)^{\frac{q^{*}}{q}\frac{p-q}{p}}\notag\\
&+C\bigg(1+\int_{\{| u_n|\leq K\}}|\nabla u_n|^{q}
\text{d}x\bigg)^{(p_1+C_1)\frac{q^*}{q}}\notag\\
&+C\bigg(1+\int_{\{| u_n|\leq K\}}|\nabla u_n|^{q}
\text{d}x\bigg)^{(p_2+C_2)\frac{q^*}{q}}
.\label{3.24}
\end{align}
Note that $p<N\Leftrightarrow\frac{q^{*}}{q}\frac{p-q}{p}<1$ and $(p_i+C_i)\frac{q^*}{q}<1 $ by \eqref{pici}. It follows from \eqref{3.24} that for $k_{0}$ large enough, $\int_{\{|u_n|\leq K\}}|\nabla u_n|^q\text{d}x $ is bounded independently of $n$ and $K$. Using \eqref{3.22} and \eqref{3.23}, we deduce that $\int_{\{|u_n|\leq K\}}|u_n|^{q^*}\text{d}x$ is also bounded independently of $n$ and $K$. Letting $K\rightarrow \infty$, we deduce that $\|\nabla u_n\|_{q}$ and $\| u_n\|_{q^{*}}$ are uniformly bounded independently of $n$. Particularly, $u_n$ is bounded in $W^{1,q}(\Omega)$. Therefore, there exists a subsequence of $\{u_n\}$ and a function $v\in W^{1,q}(\Omega)$ such that $u_n\rightharpoonup v$\ \text{weakly\ in} \ $W^{1,q}(\Omega)$,
$u_n\rightarrow v$\ strongly in $L^{q}(\Omega)$ and a.e. in $\Omega$. By Lemma~\ref{Lemma 3.2}, $u_n\rightarrow u$ in measure in $\Omega$, we conclude that $u=v$ and $u\in W^{1,q}(\Omega)$.\ $\square$

\begin{lemma}\label{Lemma 3.5}
There exists a subsequence of $\{u_n\}$ and a measurable function $u$ such that $\nabla u_n$ converges almost everywhere in $\Omega$ to $\nabla u$.
\end{lemma}
\textbf{Proof.}\ \ The proof is quite similar to Theorem 4.1 in \cite{ABFOT2003}, which can be also found in \cite{B1996}. Here we sketch only the main steps due to slight modifications. For $r_2>1$, let $\lambda=\frac{q}{pr_2}<1$, where $q$ is the same as in Lemma \ref{Lemma 3.4}.
Define $A(x,u,\xi)=\frac{a(x,\xi)}{(1+|u|)^{\theta(p-1)}}$ (for the sake of simplicity, we omit the dependence of $A(x,u,\xi)$ on $x$) and
\begin{align*}
I(n)=\int_{\Omega}((A(u_n,\nabla u_n)-A(u_n,\nabla u))\cdot\nabla (u_n-u))^{\lambda}\text{d}x.
\end{align*}
We fix $k>0 $ and split the integral in $I(n)$ on the sets $\{|u|>k\}$ and $\{|u|\leq k\}$, obtaining
\begin{align*}
I_{1}(n,k)=\int_{\{|u|>k\}}((A(u_n,\nabla u_n)-A(u_n,\nabla u))\cdot\nabla (u_n-u))^{\lambda}\text{d}x,
\end{align*}
and
\begin{align*}
I_{2}(n,k)=\int_{\{|u|\leq k\}}((A(u_n,\nabla u_n)-A(u_n,\nabla u))\cdot\nabla (u_n-u))^{\lambda}\text{d}x.
\end{align*}
For $I_{2}(n,k)$, one has
\begin{align*}
I_{2}(n,k)\!\leq\! I_{3}(n,k)\!=\!\!\int_{\Omega}\!((A_{n}(u_{n},\nabla
u_{n})-A_{n}(u_{n},\nabla
T_{k}(u)))\cdot\nabla(u_{n}-T_{k}(u)))^{\lambda}{d}x.
\end{align*}
Fix $ h>0 $ and split $ I_{3}(n,k) $ on the sets $
\{|u_{n}-T_{k}(u)| > h\} $ and $ {\{|u_{n}-T_{k}(u)|\leq h\}} $,
obtaining
\begin{align*}
I_{4}(n,k,h)=\!\!\int_{\{|u_n-T_k(u)|>h\}}\!\!\!\!\!\!\!((A_n(u_n,\nabla
u_n)-A_n(u_n,\nabla
T_k(u)))\cdot\nabla(u_n-T_k(u)))^\lambda\text{d}x,\notag
\end{align*}
and
\begin{align*}
I_{5}(n,k,h) \!&= \!\!\int_{\{|u_n-T_k(u)|\leq h\}}\!\!\!\!\!\!\!\!((A_n(u_n,\nabla
u_n)-A_n(u_n,\nabla
T_k(u)))\cdot\nabla(u_n-T_k(u)))^\lambda\text{d}x\notag
\\&=\int_\Omega((A_n(u_n,\nabla u_n)-A_n(u_n,\nabla
T_k(u)))\cdot\nabla T_h(u_n-T_k(u)))^\lambda\text{d}x\notag
\\&\leq|\Omega|^{1-\lambda}\bigg(\!\int_{\Omega}\!(A_n(u_n,\!\nabla u_n)-A_n(u_n,\!\nabla
T_k(u)))\cdot\!\nabla T_h(u_n\!-\!T_k(u))\text{d}x\!\bigg)^\lambda\notag \\&=
|\Omega|^{1-\lambda}\big(I_6(n,k,h)\big)^\lambda.\notag
\end{align*}
 For $I_{6}(n,k,h)$, it can be split as the difference $I_{7}(n,k,h)-I_{8}(n,k,h)$
\\ where
\begin{align*}
I_{7}(n,k,h)=\int_\Omega A(u_n,\nabla u_n)\cdot\nabla
T_h(u_n-T_k(u))\text{d}x,\notag
\end{align*}
\\and
\begin{align*}
I_{8}(n,k,h)=\int_\Omega A(u_n,\nabla T_k(u))\cdot\nabla
T_h(u_n-T_k(u))\text{d}x.\notag
\end{align*}
 Note that $|\nabla u_n|$ is bounded in $L^q(\Omega)$ and $\lambda p
 r_2=q$. Thanks to Lemma \ref{Lemma 3.2} and Lemma \ref{Lemma 3.4}, one may get in the same way as Theorem 4.1 in \cite{ABFOT2003} that
\begin{align*}
\lim_{k\rightarrow\infty}\!\limsup_{n\rightarrow\infty}I_1(n,k)=0,
\lim_{h\rightarrow\infty}\!\limsup_{k\rightarrow\infty}\limsup_{n\rightarrow\infty}\!I_4(n,k,h)=0,
\lim_{n\rightarrow\infty}\!\!I_8(n,k,h)=0.\notag
\end{align*}
 For $I_7(n,k,h)$, let $k> \max\{\|g\|_\infty,\|\psi\|_{\infty}\}$ and $ n\geq h+k$. Take
 $T_k(u)$ as a test function for \eqref{1.9}, obtaining
\begin{align*}
I_{7}(n,k,h)\leq\int_\Omega f_nT_h(u_n-T_k(u))\text{d}x+\int_\Omega b|u_n|^{r-2}u_nT_h(u_n-T_k(u))\text{d}x.\notag
\end{align*}
 Note that $r-1<q^*$, and $\int_{\Omega}|u_n|^{q^*}\text{d}x$ is uniformly bounded (see the proof of Lemma \ref{Lemma 3.4}),thus $|u_n|$ converges strongly in $L^{1}(\Omega)$. Therefore we have
 \begin{align*}
\lim_{n\rightarrow\infty}\int_\Omega |u_n|^{r-2}u_nT_h(u_n-T_k(u))\text{d}x=\int_\Omega |u|^{r-2}uT_h(u-T_k(u))\text{d}x.\notag
\end{align*}
 Then using the strong convergence of $f_n$ in $L^1(\Omega)$, one has
 \begin{align*}
\lim_{n\rightarrow\infty}I_7(n,k,h)
\leq\int_\Omega-fT_h(u-T_k(u))\text{d}x+\int_\Omega b|u|^{r-2}uT_h(u-T_k(u))\text{d}x.\notag
\end{align*}
It follows
\begin{align*}
\lim_{k\rightarrow\infty}\lim_{n\rightarrow\infty}I_7(n,k,h)\leq 0.\notag
\end{align*}
 Putting together all the limitations and noting that $I(n)\geq 0$, we have
\begin{align*}
\lim_{n\rightarrow\infty}I(n)=0.\notag
\end{align*}
 The same arguments as $\cite{ABFOT2003}$ give that, up to subsequence, $\nabla u_n(x)\rightarrow\nabla u(x)\
 a.e.$.\ $\square$\\

\textbf{Proof of Proposition  \ref{Proposition 3.1}.}\ \ We shall prove that $\nabla
u_n$ converges strongly to $\nabla u_n$ in $L^q(\Omega)$ for each
$q$, being given by \eqref{Q}. To do that,we
will apply Vitalli's Theorem, using the fact that by Lemma
\ref{Lemma 3.4}, $\nabla u_n$ is bounded in $L^q(\Omega)$ for each
$q$ given by \eqref{Q}. So let
$s\in(q,\frac{N(p-1)(1-\theta)}{N-1-\theta(p-1)})$ and
$E\subset\Omega$ be a measurable set. Then, we have by H\"{o}lder's inequality.
\begin{align*}
\int_E|\nabla u_n|^q\text{d}x\leq\bigg(\int_E|\nabla
u_n|^r\text{d}x\bigg)^\frac{q}{s}\cdot|E|^\frac{s-q}{s}\leq
C|E|^\frac{s-q}{s}\rightarrow 0\notag
\end{align*}
\\uniformly in $n$, as $|E|\rightarrow0$. From this and Lemma \ref{Lemma 3.5}, we
deduce that $\nabla u_n$ converges strongly to $\nabla u$ in
$L^q(\Omega)$.

 Now assume that $0\leq \theta< \min\{\frac{1}{N-p+1},\frac{N}{N-1}-\frac{1}{p-1},\frac{p-r}{p-1}\}$. Note that since $\nabla u_n$ converges
to $\nabla u$ a.e. in $\Omega$, to prove the convergence
\begin{align*}
\frac{a(x, \nabla u_n)}{(1+|u_n|)^{\theta(p-1)}} \rightarrow \frac{a(x, \nabla u)}{(1+|u|)^{\theta(p-1)}} \ \text{strongly\ in} \ (L^1(\Omega))^{N},
\end{align*}
it suffices, thanks to Vitalli¡¯s Theorem, to show that for every measurable subset $E\subset \Omega$, $\int_{E}\big|\frac{a(x, \nabla u_n)}{(1+|u_n|)^{\theta(p-1)}}\big|\text{d}x$ converges to $0$ uniformly in $n$, as $|E|\rightarrow 0$. Note that $ p-1<\frac{N(p-1)(1-\theta)}{N-1-\theta(p-1)})$ by assumptions. For any $q\in (p-1, \frac{N(p-1)(1-\theta)}{N-1-\theta(p-1)}))$, we deduce by
H\"{o}lder's inequality
\begin{align*}
\int_E\bigg|\frac{a(x, \nabla u_n)}{(1+|u_n|)^{\theta(p-1)}}\bigg|\text{d}x&\leq \beta\int_E(j+|\nabla u_n|^{p-1})\text{d}x\notag\\
&\leq \beta \|j\|_{p'}|E|^{\frac{1}{p}}+\beta\bigg(\int_E|\nabla u_n|^q\text{d}x\bigg)^{\frac{q}{p-1}}|E|^{\frac{q-p+1}{q}}\notag\\
&\rightarrow 0\ \ \text{uniformly\  in}\ n\ \text{\ as}\ |E|\rightarrow 0.\ \ \ \ \square\\
\end{align*}

\begin{lemma}\label{Lemma 3.7}
There exists a subsequence of $\{u_n\}$ such that for all $k>0$
\begin{align*}
\frac{a(x, \nabla T_k(u_n))}{(1+|T_k(u_n)|)^{\theta(p-1)}} \rightarrow \frac{a(x, \nabla T_k(u))}{(1+|T_k(u)|)^{\theta(p-1)}} \ \ \text{strongly\ in} \ (L^1(\Omega))^{N}.
\end{align*}
\end{lemma}
\textbf{Proof.}\ \ Let $k$ be a positive number. It is well known that if a sequence of measurable functions $\{u_n\}$ with uniformly boundedness in $L^s(\Omega)(s>1)$ converges in measure to $u$, then, $u_n$ converges strongly to $u$
in $L^{1}(\Omega)$. First note that the sequence $ \{\frac{a(x, \nabla T_k(u_n))}{(1+|T_k(u_n)|)^{\theta(p-1)}} \}$ is bounded in $L^{p'}(\Omega)$. Indeed, we have by \eqref{3} and Lemma \ref{Lemma 3.1},
\begin{align*}
\int_{\Omega}\bigg|\frac{a(x, \nabla T_k(u_n))}{(1+|T_k(u_n)|)^{\theta(p-1)}}\bigg|^{p'}\text{d}x &\leq
\beta \|j\|^{p'}_{p'}+\beta \int_{\Omega}\frac{|\nabla T_k(u_n)|^{p}}{(1+|T_k(u_n)|)^{\theta p}}\text{d}x\notag\\
&\leq
\beta \|j\|^{p'}_{p'}+\beta \int_{\Omega}\frac{|\nabla T_k(u_n)|^{p}}{(1+|T_k(u_n)|)^{\theta (p-1)}}\text{d}x\notag\\
&\leq C.
\end{align*}
Next, it suffices to show that there exists a subsequence of $\{u_n\}$ such that
\begin{align*}
\frac{a(x, \nabla T_k(u_n))}{(1+|T_k(u_n)|)^{\theta(p-1)}} \rightarrow \frac{a(x, \nabla T_k(u))}{(1+|T_k(u)|)^{\theta(p-1)}} \ \text{in\ measure}.
\end{align*}
Note that $u_n,u\in W^{1,q}(\Omega)$, where $q$ is the same as in Proposition \ref{Proposition 3.1}. The a.e. convergence of $u_n$ to $u$ and the fact that $ \nabla u_n\rightarrow \nabla u$ in measure imply that
\begin{align*}
\nabla T_k(u_n)\rightarrow \nabla T_k(u) \ \ \text{in\ measure}.
\end{align*}
Let $s,\ t$ be positive numbers and write $\nabla_Au=\frac{a(x, \nabla u)}{(1+|u|)^{\theta(p-1)}}$. Define
\begin{align*}
&E_{n}=\{|\nabla_A T_{k}(u_n)-\nabla_A T_{k}(u)|>s\},\\
&E_{n}^1=\{|\nabla T_{k}(u_n)|>t\},\\
&E_{n}^2=\{|\nabla T_{k}(u)|>t\},\\
&E_{n}^3=E_{n}\cap\{|\nabla T_{k}(u_n)|\leq t\}\cap\{|\nabla T_{k}(u)|\leq t\}.
\end{align*}
Note that $ E_{n}\subset E_{n}^1\cup E_{n}^2\cup E_{n}^3$ for each $n\geq 1$. Using the fact by Lemma \ref{Lemma 3.4}, the sequence $\{u_n\}$ and the fucntion $u$ are uniformly bounded in $W^{1,q}(\Omega)$, we obtain
\begin{align*}
&\mathcal {L}^{N}(E_{n}^1)\leq \frac{1}{t^q}\int_{\Omega}|\nabla T_{k}(u_n)|^{q}\text{d}x\leq \frac{1}{t^q}\int_{\Omega}|\nabla u_n|^{q}\text{d}x \leq \frac{C}{t^q},\\
&\mathcal {L}^{N}(E_{n}^2)\leq \frac{1}{t^q}\int_{\Omega}|\nabla T_{k}(u)|^{q}\text{d}x\leq \frac{1}{t^q}\int_{\Omega}|\nabla u|^{q}\text{d}x \leq \frac{C}{t^q}.
\end{align*}
 We deduce that for any $\varepsilon>0$ there exists $t_{\varepsilon}>0$ such that
 \begin{align}
\mathcal {L}^{N}(E_{n}^1)+\mathcal {L}^{N}(E_{n}^2)<\frac{\varepsilon}{3},\ \ \ \ \forall t\geq t_{\varepsilon},\ \forall n\geq 1.\label{3.39}
\end{align}
Note that for $a\geq b\geq 0,\gamma\geq 0$, we have the following inequality
 \begin{align*}
\bigg|\frac{1}{(1+a)^\tau}-\frac{1}{(1+b)^\tau}\bigg|=
\bigg|\frac{\tau(b-a)}{(1+c)^{1+\tau}}\bigg|
\leq \tau|b-a|
\ \ \ \text{for\ some}\ c\in [b,a].
\end{align*}
We deduce from \eqref{6} and \eqref{3} that in $E_{n}^3$
\begin{align*}
s&<|\nabla_A T_{k}(u_n)-\nabla_A T_{k}(u)|= \bigg|\frac{a(x,\nabla T_{k}(u_n))-a(x,\nabla T_{k}(u))}{(1+|T_{k}(u_n)|)^{\theta(p-1)}}\notag\\
&+\bigg(\frac{1}{(1+|T_{k}(u_n)|)^{\theta(p-1)}}-\frac{1}{(1+|T_{k}(u)|)^{\theta(p-1)}}\bigg)a(x,\nabla T_{k}(u))\bigg|\notag\\
&\leq \theta  (p-1)|T_{k}(u_n)- T_{k}(u)|\cdot \beta (j+|\nabla T_{k}(u)|^{p-1})\notag\\
&+ \gamma
\left\{
\begin{array}{l l}
  |\nabla T_{k}(u_n)-\nabla T_{k}(u)|^{p-1}, &  \mbox{if $1<p< 2$}\\
  |\nabla T_{k}(u_n)-\nabla T_{k}(u)|(1+|\nabla T_{k}(u_n)|+|\nabla T_{k}(u)|)^{p-2}, &  \mbox{if $p\geq 2$}
  \end{array} \right.\\
  &\leq C_{0}j|T_{k}(u_n)- T_{k}(u)|\notag\\
  &+C_{0}(1+t^{p-1}+t^{p-2})(|T_{k}(u_n)- T_{k}(u)|+|\nabla T_{k}(u_n)-\nabla T_{k}(u)|),
\end{align*}
which leads to $
E^{3}_{n}\subset F_1 \cup F_2,
$
with
\begin{align*}
 &F_1=\{ j|T_{k}(u_n)- T_{k}(u)|>\frac{s}{2C_{0}}\},\notag\\
 &F_2=\bigg\{|T_{k}(u_n)- T_{k}(u)|+|\nabla T_{k}(u_n)-\nabla T_{k}(u)|>\frac{s}{2C_{0}(1+t^{p-1}+t^{p-2})}\bigg\}.
\end{align*}
In $F_1$, we have
\begin{align*}
 \mathcal {L}^{N}(F_1)=\frac{2C_{0}}{s}\int_{F_1}\frac{s}{2C_{0}}\text{d}x< \frac{2C_{0}}{s}\int_{F_1}j|T_{k}(u_n)- T_{k}(u)|\text{d}x.
\end{align*}
By Lemma \ref{Lemma 3.2}, we deduce that there exists $n_0=n_{0}(S,C_{0},\varepsilon)$ such that
\begin{align}
 \mathcal {L}^{N}(F_1)\leq \frac{\varepsilon}{3},\ \ \forall\  n\geq n_0.\label{3.40}
\end{align}
For $F_{2}$, note that $F^{2}\subset F_3 \cup F_4,
$
with
\begin{align*}
 &F_3=\bigg\{|T_{k}(u_n)- T_{k}(u)|>\frac{s}{4C_{0}(1+t^{p-1}+t^{p-2})}\bigg\},\notag\\
 &F_4=\bigg\{|\nabla T_{k}(u_n)-\nabla T_{k}(u)|>\frac{s}{4C_{0}(1+t^{p-1}+t^{p-2})}\bigg\}.
\end{align*}
Using the convergence in measure of $\nabla T_{k}(u_n)$ to $\nabla T_{k}(u)$ and $T_{k}(u_n)$ to $T_{k}(u)$, for $t=t_{\varepsilon}$, we obtain the existence of $n_{1}=n_{1}(s,\varepsilon)\geq 1$ such that
\begin{align}
 \mathcal {L}^{N}(F_2)\leq  \mathcal {L}^{N}(F_3) + \mathcal {L}^{N}(F_4) < \frac{\varepsilon}{3},\ \ \forall\  n\geq n_1.\label{3.41}
\end{align}
Combining \eqref{3.39}, \eqref{3.40} and \eqref{3.41}, we obtain
\begin{align*}
 \mathcal {L}^{N}(\{|\nabla_{A}T_k(u_n)-\nabla_{A}T_k(u)|>s\})< \varepsilon,\ \ \forall\  n\geq \max\{n_0,n_1\}.
\end{align*}
Hence the sequence $\{\nabla_{A}T_k(u_n)\}$ converges in measure to $\nabla_{A}T_k(u)$ and the lemma follows. $\square$\\

\section{Proof of the main result}

Now we have gathered all the lemmas needed to prove the existence of an entropy solution to the obstacle problem associated with $(f,\psi,g)$. In this part, let $f_n$ be a sequence of smooth functions converging
strongly to $f$ in $L^1(\Omega)$, with $\|f_n\|_{1}\leq \|f\|_{1}+1$. We consider the sequence of approximated
obstacle problems associated with $(f_n,\psi,g)$. As mentioned in Section 2, one may proceed exactly as in \cite{DD2010}, where the existence of a solution for nonlinear elliptic equations with
degenerate coercivity was established, to get a solution $u_n\in K_{g,\psi}$ associated with $(f_n,\psi,g)$. Obviously, $u_n$ is also an entropy solution. \\

 \textbf{Proof of Theorem \ref{Theorem 1.1}.}\ \ Let $v\in K_{g,\psi}$. Taking $v$ as a test function in \eqref{1.9} associated with $(f_n,\psi,g)$, we get
\begin{align*}
\int_{\Omega}\frac{a(x, \nabla u_n)}{(1+|u_n|)^{\theta(p-1)}}\cdot\nabla (T_t(u_n-v))\text{d}x+\int_{\Omega}b|u_n|^{r-2}u_nT_t(u_n-v)\text{d}x\notag\\ \leq \int_{\Omega}f_nT_t(u_n-v)\text{d}x.
\end{align*}
Since $\{|u_n-v|<t\}\subset \{|u_n|<s\}$ with $s=t+\|v\|_{\infty}$, the previous inequality can be written
as
\begin{align}
\int_{\Omega}\!\!\chi_n\!\nabla_{\!\!A} T_{s}(u_n)\cdot\!\nabla v\text{d}x\!&\geq\!\! \int_{\Omega}\!\!-f_nT_t(u_n-v)\text{d}x+\int_{\Omega}b|u_n|^{r-2}u_nT_t(u_n-v)\text{d}x\notag\\
&+\!\!\!\int_{\Omega}\!\chi_n\!\nabla_{\!\!A} T_{s}(u_n)\cdot\!\nabla T_{s}(u_n)\text{d}x,\label{4.1}
\end{align}
where $\chi_n=\chi_{\{|u_n-v|<t\}}$ and $\nabla_Au=\frac{a(x, \nabla u)}{(1+|u|)^{\theta(p-1)}}$. It is clear that $\chi_n \rightharpoonup \chi$ weakly* in $L^{\infty}(\Omega)$. Moreover, $\chi_n$ converges a.e. to $ \chi_{\{|u-v|<t\}}$ in $\Omega\setminus\{|u-v|=t\}$. It follows that
\begin{align*}
\chi=\left\{
\begin{array}{l l}
  1, &  \mbox{in} \ \{|u-v|<t\},\\
  0, &  \mbox{in} \ \{|u-v|>t\}.
  \end{array} \right.
\end{align*}
 Note that we have $\mathcal {L}^N(\{|u-v|=t\})=0$ for a.e. $t\in (0,\infty)$. So there exists a measurable set $\mathcal {O}\subset (0,\infty)$ such that $\mathcal {L}^N(\{|u-v|=t\})=0$ for all $t\in (0,\infty)\setminus\mathcal {O}$. Assume that $t\in (0,\infty)\setminus\mathcal {O}$. Then $\chi_n$ converges weakly* in $L^{\infty}(\Omega)$ and a.e. in $\Omega$ to $  \chi=\chi_{\{|u-v|<t\}}$. Since $\nabla T_{s}(u_n)$ converges a.e. to
$\nabla T_{s}(u)$ in $\Omega$ (Proposition \ref{Proposition 3.1}), we obtain by Fatou's Lemma
 \begin{align}
\liminf_{n\rightarrow \infty}\int_{\Omega}\chi_n\nabla_A T_{s}(u_n)\cdot\nabla T_{s}(u_n)\text{d}x\geq \int_{\Omega}\chi\nabla_A T_{s}(u)\cdot\nabla T_{s}(u)\text{d}x.\label{4.2}
\end{align}
Using the strong convergence of $\nabla_A T_{s}(u_n)$ to $\nabla_A T_{s}(u)$ in $L^{1}(\Omega)$ (Lemma \ref{Lemma 3.7}) and the weak* convergence of $\chi_n$ to $\chi$ in $L^{\infty}(\Omega)$, we obtain
 \begin{align}
\lim_{n\rightarrow \infty}\int_{\Omega}\chi_n\nabla_A T_{s}(u_n)\cdot\nabla v\text{d}x= \int_{\Omega}\chi\nabla_A T_{s}(u)\cdot\nabla v\text{d}x.\label{4.3}
\end{align}
Moreover, due to the strong convergence of $f_n$ to $f$  and $|u_n|^{r-2}u_n\ $ to $|u|^{r-2}u$ (by $r-1<q^*$ and the boundedness of $\|u_n\|_{q^*}$) in $L^{1}(\Omega)$, and the weak* convergence of $T_{t}(u_n-v)$ to $T_{t}(u-v)$ in $L^{\infty}(\Omega)$, by passing to the limit in \eqref{4.1} and taking into account \eqref{4.2}-\eqref{4.3}, we obtain
  \begin{align*}
 \int_{\Omega}\chi\nabla_A T_{s}(u)\cdot\nabla v\text{d}x- \int_{\Omega}\chi\nabla_A T_{s}(u)\cdot\nabla T_{s}(u)\text{d}x&\geq \int_{\Omega}-fT_{t}(u-v)\text{d}x\notag\\
 &+\int_{\Omega}b|u|^{r-2}uT_t(u-v)\text{d}x,
\end{align*}
which can be written as
  \begin{align*}
\int_{\{|v-u|\leq t\}}\!\!\!\!\!\!\chi\nabla_A T_{s}(u)\cdot(\nabla v-\nabla u)\text{d}x\geq \!\! \int_{\Omega}\!\!-fT_{t}(u-v)\text{d}x+\!\int_{\Omega}\!\!b|u|^{r-2}uT_t(u-v)\text{d}x,
\end{align*}
or since $\chi=\chi_{\{|u-v|<t\}}$ and $\nabla (T_{t}(u-v))=\chi_{\{|u-v|<t\}}\nabla (u-v)$
  \begin{align*}
\int_{\Omega}\!\!\!\nabla_{\!\!A} u\!\cdot\!\!\nabla  T_t( u-v)\text{d}x+\!\!\int_{\Omega}\!\!b|u|^{r-2}uT_t(u-v)\text{d}x\!\leq \!\!\int_{\Omega}\!\!fT_{t}(u-v)\text{d}x,\forall\ t\in \! (0,\infty)\!\setminus\!\mathcal {O}.
\end{align*}
For $t\in \mathcal {O}$, we know that there exists a sequence $\{t_{k}\}$ of numbers in $ (0,\infty)\setminus\mathcal {O}$ such that $t_k\rightarrow t$ due to $|\mathcal {O}|=0$. Therefore, we have
 \begin{align}
\int_{\Omega}\nabla_{\!\!A} u\!\cdot\! \nabla T_{t_k}( u-v)\text{d}x+\!\int_{\Omega}\!b|u|^{r-2}uT_{t_k}(u-v)\text{d}x\leq \!\! \int_{\Omega}\!fT_{t_k}(u-v)\text{d}x.\label{4.4}
\end{align}
Since $\nabla (u-v)=0$ a.e. in $\{|u-v|=t\} $, the left-hand side of \eqref{4.4} can be written as
\begin{align*}
\int_{\Omega}\nabla_A u\cdot \nabla T_{t_k}( u-v)\text{d}x=\int_{\Omega\setminus\{|u-v|=t\}}\chi_{\{|u-v|<t_k\}}\nabla_A u\cdot \nabla( u-v)\text{d}x.
\end{align*}
The sequence $ \chi_{\{|u-v|<t_k\}}$ converges to $ \chi_{\{|u-v|<t\}}$ a.e. in $\Omega\setminus\{|u-v|=t\} $ and therefore
converges weakly* in $L^{\infty}(\Omega\setminus\{|u-v|=t\})$. We obtain
\begin{align}
\lim_{k\rightarrow \infty}\int_{\Omega}\nabla_A u\cdot \nabla T_{t_k}( u-v)\text{d}x&=\int_{\Omega\setminus\{|u-v|=t\}}\chi_{\{|u-v|<t\}}\nabla_A u\cdot \nabla( u-v)\text{d}x\notag\\
&=\int_{\Omega}\chi_{\{|u-v|<t\}}\nabla_A u\cdot \nabla( u-v)\text{d}x\notag\\
&=\int_{\Omega}\nabla_A u\cdot \nabla T_t(u-v)\text{d}x.\label{4.5}
\end{align}
For the right-hand side of \eqref{4.4}, we have
\begin{align}
\bigg|\int_{\Omega}\!\!fT_{t_{k}}(u-v)\text{d}x-\!\!\int_{\Omega}fT_{t}(u-v)\text{d}x\bigg|\leq |t_k-t|\cdot\|f\|_1\rightarrow 0\ \ \text{as}\ k\rightarrow \infty.\label{4.6}
\end{align}
Similarly, we have
\begin{align}
\bigg|\!\!\int_{\Omega}\!\!|u|^{r-2}uT_{t_k}(u-v)\text{d}x\!-\!\!\int_{\Omega}\!\!|u|^{r-2}uT_{t}(u-v)\text{d}x\bigg|\!\!&\leq \!\!|t_k-t|\!\cdot\|\!|u|^{r-1}\|_1\notag\\
&\rightarrow 0\ \text{as}\ k\rightarrow \infty.\label{4.7}
\end{align}
It follows from \eqref{4.4}-\eqref{4.7} that we have the inequality
\begin{align*}
\int_{\Omega}\!\!\nabla_{\!\!A} u\cdot \!\nabla T_{t}( u-v)\text{d}x\!+\!\!\int_{\Omega}\!\!b|u|^{r-2}uT_{t}(u-v)\text{d}x\!\leq \!\!\int_{\Omega}\!\!fT_{t}(u-v)\text{d}x,\ \forall\ t\in(0,\infty).
\end{align*}
Hence, $u$ is an entropy solution of the obstacle problem associated with $(f,\psi,g)$.

 For the regularity of an entropy solution $u$, it suffices to replace  $u_n$ and $f_n$ with $u$ and $f$ respectively in Lemma \ref{Lemma 3.4} and one can obtain the boundedness of $u$ in $W^{1,q}(\Omega)$. $\square$\\


\vspace{0.5cm}
 \noindent{\bf Acknowledgments.}  This work is supported in part by the Fundamental Research Funds for the
Central Universities: 10801B10096018.

\end{document}